\documentclass[12pt,a4paper]{amsart}

\usepackage{pdflscape}
\usepackage{amsmath}
\usepackage{amssymb}
\usepackage{latexsym}
\usepackage{srcltx}
\usepackage{graphics}
\usepackage{color}
\usepackage{epsfig}
\usepackage{color}
\usepackage{graphics}
\usepackage{graphicx}
\usepackage{epsfig}
\usepackage[ruled,vlined]{algorithm2e}

\usepackage[a4paper,hmargin=2.5cm,vmargin=2.75cm]{geometry}
\newcommand{\R}{{\mathbb R}}
\newcommand{\re}{{\mathbb R}}
\newcommand{\n}{{\mathbb N}}

\newcommand{\z}{{\mathbb Z}}
\newcommand{\cA}{{\mathcal{A}}}

\newcommand{\cT}{{\mathcal{T}}}

\newcommand{\cS}{{\mathcal{S}}}

\newcommand{\odin}{{\boldsymbol 1}}

\newcommand{\ba}{{\boldsymbol a}}
\newcommand{\bb}{{\boldsymbol b}}
\newcommand{\bc}{{\boldsymbol c}}

\newcommand{\prf}{\textbf{Proof.} }
\newcommand{\eop}{\hfill $\Box$}

\newtheorem{theorem}{Theorem}
\newtheorem{prop}{Proposition}
\newtheorem{lemma}{Lemma}
\newtheorem*{cor*}{Corollary}
\newtheorem{remark}{Remark}

\begin{document}

\date{\today}

\author{Sergei V. Konyagin}
\address{\parbox{\linewidth}{Steklov Mathematical Institute of RAS, 8 Gubkina St., 119991 Moscow, Russia}}
\email[]{konyagin@mi-ras.ru}

\author{Vladimir Yu.~Protasov}
\address{\parbox{\linewidth}{DISIM, University of L'Aquila, 67100 L'Aquila, Italy\\
Moscow State University, Moscow, 119991, Russia}}
\email[]{vladimir.protasov@univaq.it}

\author{Alexey L. Talambutsa}
\address{\parbox{\linewidth}{Steklov Mathematical Institute of RAS, 8 Gubkina St., 119991 Moscow, Russia\\
HSE University, Laboratory of Theoretical Computer Science, \\ 11 Pokrovsky Blvd., 109028 Moscow, Russia}}
\email[]{altal@mi-ras.ru}

\title[Unique expansions in number systems via refinement equations]{Unique expansions in number systems \\ via refinement equations}


\maketitle

\begin{abstract}
Using the subdivision schemes theory, we develop a criterion to check if
any natural number has at most one representation in the $n$-ary
number system with a set of non-negative integer digits $A=\{a_1, a_2,\ldots, a_n\}$ that contains zero.
This uniqueness property is shown to be equivalent to a certain restriction
on the roots of the trigonometric polynomial $\sum_{k=1}^n e^{-2\pi i a_k t}$.
From this criterion, under a natural condition of irreducibility for $A$, we deduce
that in case of prime $n$ the uniqueness holds if and only if the digits of $A$
are distinct modulo $n$, whereas for any composite $n$ we show that the latter condition is not necessary. We also establish the connection of this uniqueness to the semigroup
freeness problem for affine integer functions of equal integer slope; this together with the two mentioned criteria allows to fill the gap in the work of D.~Klarner on the question of P.~Erd\"os about densities of affine integer orbits and establish a simple algorithm to check the freeness and the positivity of density when the slope is a prime number.
\end{abstract}

\section{Number systems and semigroup freeness}
\label{sec:intro}

Given an integer $n\geq 2$, we consider a set of nonnegative integers
$A=\{a_1, \ldots , a_m\}$, which are called digits of the \emph{(non-standard) $n$-ary number system} $A$. The main results of the paper are obtained assuming $m=n$, but this equality is not necessary in general.

For an integer $k\ge 0$, we call its \emph{$A$-expansion} any tuple of digits $(\alpha_0,\ldots,\alpha_J)$ such that
\begin{equation}\label{eq.part}
k \quad = \quad \sum_{j =0}^{J} \, \alpha_{j} \, n^j,\, \text{ where } \alpha_{J} \ne 0.
\end{equation}
In particular, for $k=0$ its $A$-expansion is an empty tuple.
We say that the $n$-ary number system $A$ has \emph{uniqueness property} if any integer $k\geq 0$ has at most one $A$-expansion.

\smallskip

The uniqueness property has been studied in the literature mainly due to two applications. The first is the theory of self-similar tiles, where uniqueness gives the existence criterion for the tile. See~\cite{JT, LW, Li}
and references therein for the general multivariate case.
The second one is related to the coding theory, but as we shall show, it also turns out to be connected to the freeness of particular semigroups and densities of integer orbit sets arising from the action of integer affine functions. 

\smallskip
The active analysis of number systems with positive digits was initiated in the early 1980s by
H.A.~Maurer, A.~Salomaa and D.~Wood  in \cite{MSW}, where they introduced L-codes and discovered that if the code alphabet is unary, then an L-code corresponds to a non-standard number system\footnote{The condition $\alpha_J\ne 0$ for the leading digit in \eqref{eq.part} is set to generalize the notion of uniqueness from \cite{MSW} to the number systems containing zero, for example to the standard digit system $\{0,1,2,\ldots, n-1\}$.}. They have proved that such unary code can be uniquely decoded if and only
if the corresponding number system has uniqueness property. It was also proved
that this property necessarily fails if $m>n$, while if all $\alpha_i$
are distinct modulo $n$, then the property holds.
The authors also suggested a general decision problem of checking a number system,
whether it has the uniqueness property.
As shown by J.~Honkala in \cite{Hon1}, this problem is decidable. Some more
results concerning number systems were obtained within the next decade; for their description we refer the reader to \cite{Hon2}.

\smallskip

Even a little earlier, in \cite{K81} D.\,Klarner came to a study of another decision problem, which is closely related to the uniqueness of number system.
The motivation was coming from the question of Erd\"os on when the density of an orbit set $\langle F : S \rangle$ is positive (see \cite[p.23]{ErdGra} and \cite{Lag}).
Here the set $\langle F : S \rangle$ contains all images of a finite set of non-negative integers $S$ under an application of any number of integer-valued affine functions taken from the set
\begin{equation}
F=\{f_i=n_i x+a_i, i=1,\ldots,m \}, \text{ with $n_i\geq~2$ and $a_i\in \n\cup \{0\}$  }.
\label{eq.funcset}
\end{equation}
Klarner studied a particular case when all slopes $n_i$ are equal to some
integer $n\geq 2$ and noted that for $m=n$ the (upper) density
\[\delta(\langle F : S \rangle)=\lim\sup_{T\to \infty}\frac{\langle F : S \rangle\cap\{0,\ldots, T\}}{T+1}\]
is equal to $0$ if and only if the set $F$ has a non-trivial semigroup relation, i.e.
\begin{equation}
f_{u_1}\circ \ldots\circ f_{u_p} = f_{v_1}\circ \ldots\circ f_{v_q},
\label{eq.rel}
\end{equation}
for some $u_1,\ldots, u_p,v_1,\ldots, v_q\in \{1,\ldots,m\}$ such that $(u_1,\ldots, u_p)\ne (v_1,\ldots, v_s)$.

\smallskip

The decision problem in which, for an input set $F$, one needs to answer whether a relation of the form \eqref{eq.rel} exists,
is known as a \emph{freeness problem for a semigroup}.\footnote{This problem is known to be undecidable for integer upper-triangular $3\times 3$ matrices (see \cite{KBS} and \cite{CHK}), and it is still open for integer affine functions, which can be presented by $2\times 2$ upper-triangular matrices.}

Klarner noted, that if all slopes in \eqref{eq.funcset} are equal to $n$, the inequality $m>n$ implies that there exists a relation of type \eqref{eq.rel}; for non-equal slopes this was generalized in \cite{K82} and \cite{KT}. The important subcase $m=n$ was considered separately, but for its resolution the reader was referred to a preprint,
which has never appeared. However, it was noted that the freeness can be obtained once the set
\[\left\{0, \frac{a_2 - a_1}d, \ldots, \frac{a_m - a_1}d\right\}\]
is a complete residue system modulo $n$, where $d = \gcd(a_2 - a_1, \ldots , a_m - a_1)$.
As we will see in Section \ref{sec:res}, in the case of prime $n=m$, the freeness and this condition for residues are actually equivalent.

\smallskip

The uniqueness problem for a number system and the freeness problem for a set of integer-valued affine functions of the same slope can be related by the following observation.

\begin{prop}
Let $F=\{ f_i(x)=nx+a_i \mid  i=1,\ldots,m \}$ be the set of affine functions such that $n,a_2,\ldots,a_m \in \n$ with $n\geq 2$ and $a_1=0$.
Then the $n$-ary system $A=\{a_1,\ldots,a_m\}$ has uniqueness property if and only if the set $F$ is a free semigroup basis.
\label{p.00}
\end{prop}
\prf Suppose first that there is a semigroup relation \eqref{eq.rel}. A straightforward computation shows that
$f_{u_1}\circ \ldots\circ f_{u_p}=n^p x + c=f_{v_1}\circ \ldots\circ f_{v_q}=n^q x + d$, where
\begin{equation}
\begin{aligned}
c&=a_{u_1}+a_{u_2}n+a_{u_3}n^2+\ldots+a_{u_p}n^{p-1},\\
d&=a_{v_1}+a_{v_2}n+a_{v_3}n^2+\ldots+a_{v_q}n^{q-1}.
\end{aligned}
\label{eq.fcoeff}
\end{equation}
Since $c=d$, from \eqref{eq.fcoeff} we get two non-equal $A$-expansions of the same number.

\smallskip

Suppose that there exist two tuples $(a_{u_1},\ldots a_{u_p})$ and $(a_{v_1},\ldots a_{v_q})$ giving $A$-expansions of the same number $c=d$ in \eqref{eq.fcoeff}.
Without loss of generality, we suppose $p\geq q$ and consider $S=f^{(p-q)}_{1}=n^{p-q}x$. Then we immediately obtain, that $L=f_{u_1}\circ \ldots\circ f_{u_p}=n^px+c$ and
$R=f_{v_1}\circ \ldots\circ f_{v_q}\circ S=n^px+d$, and $L=R$ is a non-trivial semigroup relation since $a_{u_p}\ne 0$ due to the condition \eqref{eq.part}.

\eop

The restriction $a_1=0$ in the statement of Proposition \ref{p.00} is actually somewhat artificial,
and it can be bypassed by the following tool:

\begin{prop}
Let $F=\{ f_i(x)=nx+a_i \mid  i=1,\ldots,m \}$ be a set of affine functions such that $n,a_1,a_2,\ldots,a_m \in \z$ and $n\geq 2$.
Then for any $s\in \z$ the set $F$ is a free semigroup basis if and only if the set $F=\{ \tilde{f}_i(x)=nx+(a_i-s) \mid  i=1,\ldots,m \}$ is a free semigroup basis.
\label{p.01}
\end{prop}
\prf
For the function $g(x)=x-c$,  we have $g^{-1}\circ f_i \circ g=nx+(a_i-cn+c)$.
Let us take $c=s/(n-1)$, then $g^{-1}\circ f_i \circ g=\tilde{f}_i$. Consequently, the existence of a group relation
\[
f_{u_1}\circ \ldots\circ f_{u_p} = f_{v_1}\circ \ldots\circ f_{v_q}.
\]
is equivalent to the existence of the conjugated relation
\[
g^{-1}\circ f_{u_1}\circ \ldots\circ f_{u_p} \circ g = g^{-1}\circ f_{v_1}\circ \ldots\circ f_{v_q} \circ g,
\]
which can be also written as
\[
\tilde{f}_{u_1}\circ \ldots\circ \tilde{f}_{u_p} = \tilde{f}_{v_1}\circ \ldots\circ \tilde{f}_{v_q}.
\]
\eop

\medskip

Propositions \ref{p.00} and \ref{p.01} together with Theorem \ref{th.20} from Section \ref{sec:unique_crit} will provide an algorithm to check freeness for the set $F$ in the case $m=n$, which was missing in Klarner's considerations. If $n$ is prime, the freeness check for the set $F$ (and for the density $\delta(\langle F : S \rangle$) to be positive) can be done much easier by the use of Theorem~\ref{th.prime} from Section \ref{sec:res}.

\section{Subdivision scheme and transition operator}
\label{sec:refin_eq}

Subdivision schemes are iterative algorithms for the linear approximation of functions
from their values on a mesh. They originated in the late 1980s with N.~Dyn, D.~Levin, J.~Gregory, C.~De Boor, S.~Dubuc, etc., see~\cite{CDM, DGL, DL} for many references.
They are widely used for interpolation and approximation of smooth functions and
in modeling of curves and surfaces. Their applications to combinatorics and number theory are also known, see \cite{FS, LW, P00p, P04b}. In paper \cite{P17} the refinement schemes were applied to the problem of representing numbers in the binary number system with an arbitrary set of digits.

Some special cases of subdivisions appeared earlier in the works of
G.~de~Rham (cutting angle scheme)~\cite{dR1,dR2} G.M.~Chaikin~\cite{Cha}, C.A.~Micchelli and H.~Prautzsch \cite{MP} etc. 

We consider only stationary univariate schemes defined by
an integer  contraction factor~$n\ge 2$ and by a finitely supported sequence of real numbers $\{c_i\}_{i\in \z}$.
After  a suitable shift of numeration it can be assumed that~$c_i = 0$
for all~$i\notin \{0, \ldots , N\}$ and~$c_0c_N \ne 0$.
We denote by~$l_{\infty}$ the space of bounded sequences
$\{x_k\}_{k\in \z}$ and by~$\delta \in l_{\infty}$ the sequence~$(\delta )_k  = \delta_{0k}$
(the Kronecker symbol).

\smallskip

The {\em subdivision operator} $S: l_{\,\infty} \to l_{\,\infty}$ acts
on~$l_{\,\infty}$ as follows:
\begin{equation}\label{eq.subd}
(S g)_k \quad = \quad \sum_{i \in \z}\, c_{k-ni}\, g_i
\end{equation}
where $g = (g_i)_{i \in \mathbb Z} \in l_{\infty }$.  Furthermore, we consider the {\em transition operator}~$T$ on the space of compactly supported tempered distributions~$\cS_0'$:
\begin{equation}\label{eq.trans}
\bigl[Tf\bigr](t)\quad = \quad \sum_{k\in \z} \, c_k \, f(n t-k)\, , \qquad f\, \in \, \cS_0'\,.
\end{equation}
There is a simple relation between the subdivision scheme and the transition operator.
For every~$j\in \n$, we have:
\begin{equation}\label{eq.subd-trans}
[T^jf](t)\quad = \quad \sum_{k\in \z} \ \bigl( \, S^j\delta \, \bigr)_k \ f(n^j \, t\, -\, k)\, , \qquad f\in \cS'\, .
\end{equation}
The proof can be found, for example, in~\cite{CDM} or verified by a direct computation.
The following theorem is well known~\cite{NPS}:
\smallskip

\noindent \textbf{Theorem A.} {\em Assume~$\sum_k c_k = n$. Then for every $f \in \cS_0'$ such that~$(f, \odin)  = 1$, the
 sequence~$T^jf$ converges in~$\cS'$ to a unique solution~$\varphi \in \cS_0'$
 of the functional equation
 \begin{equation}\label{eq.ref}
\varphi(t) \quad = \quad \sum_{k\in \z} \, c_k\, \varphi (nt - k)\,
\end{equation}
such that~$(\varphi, \odin) = 1$
This solution is supported on the segment~$[0,N]$ and satisfies the
equation~$T\varphi = \varphi$.
 }
\smallskip

Of course, we could have omitted the normalization condition~$(f, \odin)  = 1$ and obtain by
homogeneity: for every~$f\in \cS_0'$ the sequence~$T^jf$ converges to~$(f, \odin)\, \varphi$.
The function~$\varphi$ is called in the literature~{\em refinable function} and
equation~(\ref{eq.ref}) is called a {\em refinement equation}.

To every refinement equation we associate the trigonometric
polynomial
$$
\bc(\xi) \quad = \quad \frac1n \ \sum_{k \in \z} \, c_k \, e^{-\ 2\pi i k \xi}
$$
called {\em mask}. This is the characteristic function of the sequence~$\{c_k\}_{k\in \z}$.
Computing the Fourier transform of both parts of the refinement equation gives the following equation in the frequency domain:
\begin{equation}\label{eq.ref-fourier}
\widehat \varphi(\xi) \quad = \quad \,  \bc(\xi/n)  \, \widehat \varphi(\xi/n)
\end{equation}
Iterating~$j$ times we get
\begin{equation}\label{eq.ref-fourier-j}
\widehat \varphi(\xi) \quad = \quad \,  \bc \, \bigl(n^{-1}\xi\bigr)
\cdots \bc \, \bigl(n^{-j}\xi\bigr)\, \widehat \varphi(n^{-j}\xi).
\end{equation}
Now we focus on the case of  nonnegative coefficients~$c_k$.
Some of the results below are known~\cite{P00p}, we give their proofs for the convenience of the reader. The following simple observation has been made independently in a number of
papers~\cite{D, DGL, P00p, Z}.
\begin{prop}\label{p.10}
If all~$c_k$ are nonnegative, then~$\varphi$ is a Borel measure
of pure type, namely, it is either absolutely continuous (i.e., $\varphi \in L_1$) or purely singular.
\end{prop}
\prf Choosing~$f\, = \, \chi_{[0,1]}$, we see that the
functions~$T^jf$ are all non-negative and, hence, so is their limit.
Thus, $\varphi$ is a nonnegative distribution, that is a Borel measure.

According to the Lebesgue theorem, there exists a unique
representation of~$\varphi$ in the sum of absolutely continuous
and singular measures~$\varphi = \varphi_{cont} + \varphi_{sing}$. The uniqueness implies that
both~$\varphi_{cont}$ and~$\varphi_{sing}$  satisfy the refinement equation~(\ref{eq.ref}).
However, this equation possesses a unique solution up to normalization.
Therefore, one of those functions is zero.

\eop
\medskip

Proposition~\ref{p.10} rises the question of separating the cases
of absolute continuity and singularity of the measure~$\varphi$.
The following criterion is proved by applying~(\ref{eq.ref-fourier-j}) and the Poisson summation formula:
\begin{prop}\label{p.20}
Suppose all~$c_k$ are nonnegative; then
$\varphi \in L_1$ if and only if~$\widehat \varphi (m) = 0$
for all integer~${m\ne 0}$. Moreover, in this case~$\varphi(t) \le 1$ almost everywhere.
\end{prop}
\prf (Necessity). Assume that there exists~$m\in \z\setminus \{0\}$ such
that~$\widehat \varphi (m) \ne 0$.
Taking an arbitrary number~$j\in \n$ and substituting~$\xi = n^jm$ to equation~(\ref{eq.ref-fourier-j}), we get
$$
\widehat \varphi(n^jm) \ = \ \bc  \bigl(n^{j-1}m\bigr)\cdots \bc  \bigl(m\bigr)
\widehat \varphi\bigl(m\bigr) \ = \ \widehat \varphi \bigl(m\bigr).
$$
The latter equality holds because all the numbers~$n^{j-s}m, \, s=1, \ldots , j$,  are
integer and therefore, $\, \bc \bigl( n^{j-s}m \bigr) \, = \, \bc \bigl(0 \bigr) \, = \, 1$.
Thus,~$\widehat \varphi(n^jm) \, = \, \widehat \varphi \bigl(m\bigr)$ for all~$j\in \n$.
On the other hand, if~$\varphi \in L_1(\re)$, then~$\widehat {\varphi}(\xi) \to 0$ as~$\xi\to \infty$.
This is not true for~$\xi = n^jm\, , \, j\to \infty$, hence~$\varphi \notin L_1$.
\smallskip

(Sufficiency). If~$\widehat \varphi(m) = \delta_{\, m 0} $, then
applying the Poisson summation formula to the function~$\varphi (t-\cdot)$ we get
\begin{equation}\label{eq.poiss}
\sum_{k\in \z} \varphi (t - k)\quad = \quad \sum_{m\in \z} \widehat \varphi (m)
e^{\, 2\pi i m t}
\quad = \quad \widehat \varphi (0) \ = \  1 \, .
\end{equation}
Thus, $\sum_{k\in \z} \varphi (t - k) = \odin$. Since~$\varphi \ge 0$, we see that
$\varphi \le \odin$, i.e., for every nonnegative test function~$f\in \cS$, we have
$(\varphi, f)\, \le \, \int_{\re} f\, dt$. Hence, $\varphi$ is majorized by the Lebesgue measure and, therefore $\varphi \in L_1$.

\eop
\bigskip

The criterion of Proposition~\ref{p.20} has a disadvantage that
it involves the function~$\widehat \varphi$, which is a priori unknown.
A criterion in terms of the coefficients~$\{c_k\}_{k\in \z}$ exploits the $n$-ary tree.
We define the tree~$\cT$ as follows. The root is associated to zero and has~$n-1$ children with the numbers~$\frac{k}{n}, \, k = 1, \ldots , n-1 $.
The further construction is by induction: every vertex~$\alpha$, apart from the
root, has~$n$ children~$\frac{\alpha + k}{n}, \, k = 0, \ldots , n-1 $.
The root has level zero, its children form the first level, etc. Thus,
the~$j$th root consists of $(n-1)n^{j-1}$ numbers~$\, n^{-j}k\ , k = 1, \ldots , n^{j}-1, \
k\not \equiv  0\, ({\rm mod}\, n)$.  A subset~$\cA$ of vertices of~$\cT$
is called a {\em minimal cut set} if every infinite path along the tree from the root (all paths are without backtracking, the root is not in the path) has exactly one
common vertex with~$\cA$. All minimal cut sets are finite. The
simplest one is the set of vertices of the first level: $\cA = \{\frac{1}{n}, \ldots , \frac{n-1}{n}\}$.

\medskip

\begin{prop}\label{p.30}
Suppose all~$c_k$ are nonnegative; then a function
$\varphi \in L_1$ if and only if
there exists a minimal cut set of the~$n$-ary tree~$\cT$ that consists of
roots of the mask~$\bc$.
\end{prop}
\prf
(Necessity). Since~$\varphi$ is compactly supported,
it follows from the Paley-Wiener theorem that $\widehat \varphi$ is an entire function.
Therefore, it has finitely many zeros, if any, on the unit disc.
This implies that there exists~$j\ge 1$ such that
$\widehat \varphi $ does not have zeros on the~$j$th level of~$\cT$,
i.e., $\bc(n^{-j}m) \ne 0$ for every natural $m < n^j, \, m \, \not \equiv 0 \, ({\rm mod}\, n)$.
Let~$m = \sum_{k=0}^j d_kn^k \, = \, d_j \ldots d_0$ be the standard $n$-adic expansion of~$m$, possibly, starting with zeros, $d_0 \ne 0$. Add the digit~$d_{j+1} = 0$ and
substitute~$\xi = m$ to equation~(\ref{eq.ref-fourier-j}). We obtain
\begin{equation}\label{eq.ref-fourier-m}
\widehat \varphi(m) \ = \  \,
\widehat \varphi(n^{-j}m)\, \prod_{q=0}^j \bc \, \bigl(d_{j+1}\ldots d_{j-q+1}.d_{j-q}\ldots d_0 \bigr)   \ = \
\widehat \varphi(n^{-j}m)\,
\prod_{q=0}^j \bc \, \bigl(0.d_{j-q}\ldots d_0 \bigr)  .
\end{equation}
 Since~$\widehat \varphi(m) = 0$ and~$\widehat \varphi(n^{-j}m) \ne 0$, we see that
 one of the numbers~$n^{q-j}m \, = \, 0.d_{j-q}\ldots d_0$ must be a root of~$\bc$.
 Those numbers form a finite path of length~$j$. Thus, every
 path of length~$j$ contains a root of the polynomial~$\bc$.
 \smallskip

(Sufficiency). Assume there exists a minimal cut set~$\cA \subset \cT$
such that $\bc(\cA) = 0$. Then taking an arbitrary natural number~$m$ and
applying formula~(\ref{eq.ref-fourier-m}) we obtain~$\widehat \varphi(m) =0$,
which completes the proof for positive~$m$. The proof for negative~$m$ is the same.

\eop

\section{Criterion of uniqueness for a number system}
\label{sec:unique_crit}

Now we consider a set of digits $A = \{a_1, \ldots , a_n\}$ for the (non-standard) $n$-ary number system.
We will assume that $0 = a_1 < \ldots < a_n$. For any integer~$k \ge 0$ we denote by $b(k)$
the total number of its $A$-expansions having form \eqref{eq.part} and set formally~$b(k) =0$
for any integer $k<0$.

\smallskip

The uniqueness expansion property means that~$b(k) \le 1$ for all~$k \in \z$.
We are going to see that the sequence~$\{b(k)\}_{k\in \z}$
is generated by a subdivision operator~$S$ with the following coefficients:
\begin{equation}\label{eq.c}
c_i \ = \
\left\{
\begin{array}{ccl}
1 & , & i \in A\\
0& , &  {\rm otherwise}
\end{array}
\right.
\end{equation}
This sequence will be referred to as an {\em indicator sequence of~$A$}.
We start with the following simple observation:
\begin{lemma}\label{l.10}
The sequence~$b(k)$ satisfies the following recurrent relations:
\begin{equation}\label{eq.b-rec}
b\,(nq + d) \ = \ \sum_{\,s \in \z} \, c_{\,ns + d}\ b\, (q - s)\, ,  \quad d=0, \ldots , n-1,
\qquad q\ge 0.
\end{equation}
\end{lemma}

\prf In expansion~(\ref{eq.part}) for~$k=nq+d$,  the digit~$\alpha_0 \in A$  must be
such that~$\alpha_0 \equiv d \  ({\rm mod}\,  n)$, hence $\alpha_0 \, = \, ns + d$
for some~$s$. Subtracting~$d$ from both sides of~(\ref{eq.part})   and
dividing by~$n$, we get
$q-s = \sum_{j} \alpha_j\, n^{j-1}$, provided that $ns + d \in A$, i.e., $\alpha_{ns + d} = 1$.
Thus, the total number of expansions of the number $nq + d$ is equal to the sum
of numbers of expansions of  $q-s$ over all $s$ such that $\alpha_{ns + d} =1$.
This completes the proof.

\eop
\medskip

If we denote~$k=nq+d$ and~$i=q-s$, then the equation~(\ref{eq.b-rec})
becomes
\begin{equation}\label{eq.b-rec-k}
b(k) \ = \ \sum_{\,i \in \z} \, c_{\,k - ni}\, b(i).
\end{equation}
The right hand side is precisely the subdivision operator~(\ref{eq.subd}).
Thus, we obtain
\begin{cor*}
Let~$\bb^{(j)} \in l_{\infty}$ be the sequence~$\{b(k)\}_{k=0}^{n^j-1}$, complemented by zeros
for~$k < 1$ and for $k \ge n^j$. Then~$\bb^{(j)} \, = \, S\bb^{(j-1)}$.
\end{cor*}

\begin{theorem}\label{th.10}
For every~$j\ge 1$, we have
\begin{equation}\label{eq.b-subd}
b(k)\quad = \quad  (S^j \delta)_k \ , \qquad k = 0, \ldots , n^{j}-1
\end{equation}
where~$S$ is the subdivision operator defined by the
indicator sequence of~$A$.
\end{theorem}
Thus, the sequence~$\bb^{(j)}$ defined in the Corollary coincides with the sequence~$S^j\delta$.
\smallskip

\prf We argue by induction in~$j$.
For~$j=0$, we have~$b(0) = 1$ and   $(S^0\delta)_0 = (\delta)_0 = 1$, so, this
is true. The transfer $j-1\, \mapsto \, j$ is
provided by the Corollary.

\eop
\medskip

Now we formulate the fundamental result.
Let  $\{c_k\}_{k\in \z}$ be the indicator sequence of~$A$, $\bc(\cdot)$ be the mask of this sequence, $\cT$ be the $n$-ary tree.
\begin{theorem}\label{th.20}
The following assertions are equivalent
\smallskip

\noindent a) Every natural number possesses at most one
$n$-ary expansion with the digit set~$A$;

\smallskip

\noindent b) The  compactly supported solution of the refinement
equation~$\varphi(t) = \sum_{i=1}^n \varphi(nt - a_i)$ belongs to~$L_1(\re)$;

\smallskip

\noindent c) The $n$-ary tree~$\cT$ possesses a minimal cut set
that consists of roots of the polynomial~$\bc(t) = \frac1n\, \sum_{k=1}^n e^{-2\pi i a_k t}$.

\end{theorem}
\prf The equivalence of b) and c) follows from Proposition~\ref{p.30}.
Let us prove the equivalence of~a) and~b).
\smallskip

\noindent $\mathbf{a) \, \Rightarrow \, b)}$.
Assume the contrary: the property a) holds, i.e., $b(k) \le 1, \, k\in \z$,
but~$\varphi \notin L_1$. The latter means that $\varphi$ is a purely singular
Borel measure
(Proposition~\ref{p.10}). Consider the
transition operator~$[Tf](t)\, = \, \sum_{i=1}^n \varphi(nt - a_i)$
and take $f=\chi_{[0,1]}$. By equation~(\ref{eq.subd-trans}), we have
$$
[\,T^jf\,]\,(t)\quad = \quad \sum_{k\in \z} \, \bigl( S^j\delta \, \bigr)_k \ f(n^j t\, -\, k)
\quad = \quad \sum_{k\in \z} \, \bigl( \, S^j\delta \, \bigr)_k \, \chi_{[n^{-j}k, n^{-j}(k+1)]}\, .
$$
On the other hand, $T^jf \, \to \, \varphi$ due to Theorem~A. The convergence is in~$\cS'$,
which,  for  nonnegative distributions,  means the convergence in  measure.
Applying now Theorem~\ref{th.10} we obtain~$\bigl( S^j\delta \bigr)_k = b(k)$
for all~$k \le n^j-1$. Since the segments~$[n^{-j}k, n^{-j}(k+1)]$
cover the segment~$[0,1]$ when~$k$ runs from zero to $n^j-1$, we
conclude that the restriction of the function
$$
\sum_{k\in \z} \bigl( S^j\delta \bigr)_k \, \chi_{[n^{-j}k, n^{-j}(k+1)]}
$$
to the segment~$[0,1]$ is equal to
$$
f_j\ = \ \sum_{k=0}^{n^j-1} b(k) \, \chi_{[n^{-j}k, n^{-j}(k+1)]}\, .
$$
Therefore,
\begin{equation}\label{eq.b-limit}
f_j \ = \
T^jf \, \bigl|_{[0,1]} \ \to \ \varphi \bigl|_{[0,1]}\qquad \mbox{{\rm as}}\quad j\to \infty\, .
\end{equation}
By Proposition~\ref{p.10}, the function~$\varphi$ is of pure type,
hence, the assumption~$\varphi \notin L_1$ implies that~$\varphi$
is purely singular. So, its restriction to the segment~$[0,1]$ is purely singular
as well. On the other hand, if~$b(x)\le 1$ for all~$b(k)$, then
the function on the right hand side of~(\ref{eq.b-limit})
does not exceed one and hence, so its limit~$\varphi\bigl|_{[0,1]}$. Thus,
the function~$\varphi\bigl|_{[0,1]}$ is majorized by the Lebesgue measure
and is not purely singular.
\smallskip

\noindent $\mathbf{b) \, \Rightarrow \, a)}$. Assume that there is a number that has
at least two $A$-expansions, the longest of which contains~$r$ digits. Let~$A_r \, = \, \bigl\{\sum_{j=0}^{r-1} \alpha_jn^j, \ \alpha_j\in A, \, j \le r\bigr\}$. Since $a_1=0$, the set $A_r$ contains all $A$-presentable numbers of length at most $r$.
Hence, at least two elements of~$A_r$ coincide,  so we have~$|A_r| < n^r$.
Consider the~$r$th power of the transition operator:
\begin{equation}\label{eq.trans-r}
[T^rf](t) \quad = \quad \sum_{k\in A_r}^n b_r(k) \varphi(n^r t - k)\, ,
\end{equation}
where~$b_r(k)$ is the total number of~$A$-expansions of the number~$k$
with at most~$r$ digits.
Clearly,~$T^r\varphi=\varphi$, hence,
the refinement equation with the transition operator~$T^r$
possesses the same solution~$\varphi$.
By the Hutchinson theorem~\cite{H}, there is a unique
compact set~$Q\subset \re$ such that
$Q \, = \, \bigcup\limits_{k\in A_r}n^{-r}(Q+k)$. This is a fractal generated
by affine contractions~$n^{-r}(\cdot +k)$. The Lebesgue measure
$\mu(Q)$ does not exceed~$\sum_{k\in A_r}n^{-r}\mu(Q+k)\, = \, |A_r|\, n^{-r}\, \mu(Q)$.
Since $n^{-r}|A_r| < 1$, it follows that~$\mu(Q) = 0$.
On the other hand, $T^r$ respects the set of Borel measures supported on~$Q$.
Taking an arbitrary such measure~$\psi$ normalized by the condition~$(\psi, \odin) = 1$
and applying Theorem~A, we obtain~$T^{rj}\psi \, \to \, \varphi$
as $\ j\to \infty$.
Therefore,~$\varphi$ is also supported on~$Q$ and hence, is purely singular.

\eop

\begin{remark} \normalfont
Now we see that Theorem \ref{th.20} together with Proposition \ref{p.00} and Proposition \ref{p.01} provide an algorithm to check freeness for a set $F=\{f_i=n x+a_i, i=1,\ldots,n \}$, with $n,a_1,\ldots,a_n \in \z$ and $n\geq 2$. Indeed, the degree of the trigonometric polynomial $\sum_{k=1}^n e^{-2\pi i a_k t}$ is bounded from above, so is the cardinality of its roots on the segment $[0,1]$. Thus, the possible minimal cut sets of the $n$-ary tree in Theorem \ref{th.20} also have bounded size, which in turn bounds the height of possible tree vertices in the minimal cut. It follows that the check of all such vertices can be done through their finite enumeration and polynomial evaluation in the corresponding points.
\end{remark}

\section{Uniqueness of the expansion and digits residues}
\label{sec:res}
In this Section we will use the criterion from Theorem \ref{th.20}
to study the uniqueness property for the $n$-ary number system
consisting of exactly $n$ digits. Let again $A=\{a_1, \dots, a_n\}$
be the set of nonnegative digits such that $0\leq a_1 < \dots < a_n$.
We say that $n$-ary number system $A$ is \emph{irreducible} if
\begin{equation}
\gcd(a_2 - a_1, \dots, a_n - a_1, n) = 1.
\label{tag1}
\end{equation}

Since zero is allowed to be present in the set of digits, we remind that \emph{A-expansion} of a number $k$ has the form
\begin{equation}
k = \sum_{j=0} ^J \alpha_j n^j,\, \text{ where } a_J\neq 0.
\label{tag2}
\end{equation}

If $k$ has at least one such representation, we will say that $k$
is a \emph{representable number}.

\smallskip

We will say that $n$-ary number system $A$ has \emph{weak uniqueness property} if two $A$-expansions of the same length represent different numbers. The following statement is useful for working with sets of digits that contain only positive numbers.

\begin{prop}
Let $A=\{a_1,\ldots,a_m\}$ be a set of digits such that $0<a_1<a_2<\ldots<a_m$. Then the $n$-ary system $A$ has weak uniqueness property if and only if the $n$-ary system $B=\{0, a_2-a_1,\ldots, a_m-a_1\}$ has uniqueness property.
\label{prop:wrup}
\end{prop}

\prf
We prove both directions by contraposition. First, if $A$ does not have the weak uniqueness property, then there exist two tuples $U=(u_0,u_1,\ldots,u_s)$ and $V=(v_0,v_1,\ldots,v_s)$ from $A^s$ such that
\[
\sum_{j=0}^s u_j n^j = \sum_{j=0}^s v_j n^j,
\] 
and by subtraction of $a_0(1+n+\ldots +a^s)$ from both sides we obtain 
\[
\sum_{j=0}^s (u_j-a_0) n^j = \sum_{j=0}^s (v_j-a_0) n^j,
\] 
which gives two distinct tuples $U'$ and $V'$ from $B^p$. By removing most significant digits until they are non-zero, we obtain two tuples refuting the uniqueness property for the system $B$.

\smallskip

Now, suppose that the $n$-ary system $B$ does not have the uniqueness property and there exists $U=(u_0,u_1,\ldots,u_s)\in B^s$ and $V=(v_0,v_1,\ldots,v_t)\in B^t$ such that
\begin{equation}
\sum_{j=0}^s u_j n^j = \sum_{j=0}^t v_j n^j.
\label{eq.nele}
\end{equation} 

Without loss of generality, we may suppose that $s\geq t$. Consider two tuples $U',V'\in A^{s+1}$ obtained as $U'=(u_0+a_0,u_1+a_0,\ldots,u_s+a_0)$ and $(v_0+a_0,v_1+a_0,\ldots,v_t+a_0,a_0,\ldots, a_0)$. They represent the same number in the system $A$, obtained from \eqref{eq.nele} by adding $a_0(1+n+\ldots +a^s)$.

\eop

Now we are ready to prove two results on how the uniqueness property of the $n$-ary system with $n$ digits is related to their residues modulo $n$.

\begin{theorem}\label{th.composite}
Let $n$ be a composite number. Then there exists an irreducible set $A$ consisting of $n$ non-negative digits such that
the $n$-ary number system $A$ possesses the uniqueness property, but $A$ contains two numbers congruent modulo $n$.
\end{theorem}
\prf
Since $n$ is a composite number, there exist integers $n_1,n_2$, both larger than $1$ such that $n = n_1 n_2$.
We construct the following set of digits
\begin{equation}
 A=\{un_1n + v:\,0\le u < n_2,\,0\le v < n_1\}.
 \label{eq.digset}
\end{equation}
As $n_1>1$, then by taking $u=0$ we get digits $a_1 = 0$, $a_2 = 1$, so the condition \eqref{tag1} holds.
Since $n_1$ divides $un_1n$ and $0\leq v< n_1$ in \eqref{eq.digset}, no digit can be congruent to $n_1$ modulo $n$. Hence, the digits do not form the complete system
of residues modulo $n$.

From \eqref{eq.digset} it follows that no digits $a,a'$ can satisfy the congruence $a - a'\equiv un_1 \pmod n$
for $1\le u < n_2$. Then, considering the difference of two expansions of  type \eqref{tag2} modulo $n$ we obtain that if $y,y'$ are representable numbers,
we also have
\begin{equation}
 y - y'\not\equiv un_1 \;(\bmod\;n),\,\text{ when }1\le u < n_2.
 \label{tag3}
\end{equation}

Now we will prove that any positive integer has at most one expansion of type \eqref{tag2}.
Assume the contrary, and let $x$ be the least positive integer having at least
two expansions
$$ x =  \sum_{j=0} ^J \alpha_j n^j = \sum_{j=0} ^{J'} \alpha'_j n^j\,\text{ with } \alpha_J\neq 0 \text{ and } \alpha'_{J'}\neq 0.$$
Then the last digits $\alpha_0$ and $\alpha'_0$ must be distinct: otherwise
the number $(x - \alpha_0)/n$ also has two distinct expansions, but this contradicts
the supposition of minimality for chosen $x$. Hence, $\alpha_0\neq \alpha'_0$. Without loss of generality we will presume that
$\alpha_0 < \alpha'_0$. Then we have
\[ x = ny + \alpha_0 = ny' + \alpha'_0,\]
and $y > y'$ are some representable numbers. Thus, $\alpha'_0 - \alpha_0 = n(y - y')$, hence 
\begin{equation}
\alpha'_0 - \alpha_0 \equiv 0\pmod n.
\label{eq.modeq}
\end{equation}
Then, if the digits $\alpha'_0$ and $\alpha_0$ were formed in $A$ as $\alpha'_0=u_1n_1n+v_1$ and $\alpha'=u_2n_1n+v_2$, from \eqref{eq.modeq} we obtain $v_1=v_2$, hence $\alpha'_0-\alpha_0 = (u_2-u_1) n_1n$. Therefore, $y - y' = (u_2-u_1)n_1$ with $ 1 \le u_2-u_1 < n_2$, which contradicts \eqref{tag3}.
This proves the uniqueness of the expansion.

\eop

\smallskip

\begin{theorem}\label{th.prime}
For an arbitrary prime number $p$ the following holds: an irreducible $p$-ary number system $A$ has a uniqueness property, if and only if the set $A$ does not contain two numbers which are congruent modulo $p$.
\end{theorem}
\prf
The sufficiency is straightforward (e.g., \cite[Theorem 4]{MSW}), so we prove the necessity. 

Let $a_1 < \dots < a_p$ be an irreducible set of non-negative digits such that
every number has at most one representation in base $p$ with digits $a_1, \dots, a_p $.
We will show that in this case $a_1, \dots, a_p $ form a complete
system of residues modulo $p$. To do so, we will use Theorem \ref{th.20}.
For any tuple $\bc = (c_1,\dots,c_p)\in\z_{\geq 0}^p$ we define the polynomial
$$P_\bc (z) = \sum_{j=1}^p z^{c_j}.$$
For $m\in\n$, a complex number $z$ is a primitive $p^m$th root of unity
if $z^{p^m} = 1$, but $z^{p^{m-1}}\neq 1$. Denote $\ba = (a_1,\dots,a_p)$ and $\bb=(0,a_2-a_1,\ldots,a_p-a_1)$. Using Proposition \ref{prop:wrup} we can obtain that the uniqueness property for the $n$-ary system $A$ implies the uniqueness for the system $B=\{0,a_2-a_1,\ldots,a_p-a_1\}$, to which we can apply Theorem \ref{th.20}. The direction (a)$\Rightarrow$(c) in Theorem \ref{th.20} implies that some root of the polynomial $P_\bb$ is a primitive $p^m$th root of unity for some $m\in\n$, from which it follows that for some $m$ all primitive $p^m$th roots of unity are roots of $P_\bb$. Since $P_\ba=e^{-2\pi i a_0}P_{\bb}$, the latter is also true for the polynomial $P_\ba$. This will suffice for our purpose.

Recall that $z$ is a $p^m$th primitive root of unity if and only if
it is a root of the cyclotomic polynomial
\[Q(z) = Q_{p^m}(z) = \sum_{j=0}^{p-1} z^{jp^{m-1}}.\]
Next, we observe that for any $a,\tilde a\in\z_+$ such that $a\ge\tilde a$ and $a\equiv\tilde a(\bmod p^m)$,
the polynomial $z^a - z^{\tilde a}$ is divisible by $z^{p^m}-1$ which in turn is divisible by $Q$.
Hence, all primitive $p^m$th roots of unity are roots of $z^a - z^{\tilde a}$. For any
$j=1,\dots,p$ we choose $\tilde a_j$ such that $0 \le \tilde a_j < p^m$ and $\tilde a_j\equiv a_j(\bmod p^m)$.
Denote
$$ \tilde P(z) = \sum_{j=1}^p z^{\tilde a_j}.$$
Then the primitive $p^m$th roots of unity are roots of $\tilde P$. This implies that
$ \tilde P = QR,$ where $R$ is an integral polynomial of degree $<p^{m-1}$. Moreover,
all coefficients of $R$ are nonnegative. We have $\tilde P(1) = Q(1) = p$. Hence, $R(1)  =1$,
and $R(z) = z^u$ for some $u$ with $0\le u <p^{m-1}$.

If $m=1$, then $R(z) = 1$, $\tilde P(z) = \sum_{j=0}^{p-1} z^j$. We see that the multiset
$\{\tilde a_1,\dots,\tilde a_p\}$ is actually the set $\{0,\dots,p-1\}$. This means that
$a_1, \dots, a_p $ form the complete system of residues modulo $p$ as desired.

Let $m>1$. Then all numbers $\tilde a_j$ are congruent to $u$ modulo $p^{m-1}$.
Consequently, all numbers $a_j$ are also congruent to $u$ modulo $p^{m-1}$. This does
not agree with the irreducibility of the set $\{a_1, \dots, a_p\} $.

\eop

\begin{remark} \normalfont
If $n$ is prime, Theorem~\ref{th.prime} provides an easy freeness test for a set of functions $F=\{f_i=n x+a_i, i=1,\ldots,m \}$, with $n\geq~2$ and $a_i\in \n\cup \{0\}$.  First, we increasingly order the set $A=\{a_1,\ldots,a_n\}$ and compute $d=\gcd(a_2-a_1, \ldots, a_p-a_1)$. Then, the set $\{0,(a_2-a_1)/d,\ldots,(a_p-a_1)/d\}$ is an irreducible $p$-ary number system, and we are left to check whether its elements are all distinct modulo $p$. Here, the correctness of scaling the coefficients by $d$ can be proved similarly to the correctness of the shift in Proposition~\ref{p.01} by the use of the conjugation function $g(x)=dx$.
\end{remark}

\section{Acknowledgements}
The work of V.Yu.~Protasov was prepared with the support of the Theoretical Physics and Mathematics Advancement Foundation ``BASIS'' grant No. 22-7-1-20-1. The work of A.L.~Talambutsa was prepared within the framework of the HSE University Basic Research Program and with the support of the Theoretical Physics and Mathematics Advancement Foundation ``BASIS'' grant No. 22-7-2-32-1 . 


\begin{thebibliography}{}

\bibitem{CHK}
J.\,Cassaigne, T.\,Harju, J.\,Karhum\"aki,
\newblock {\em On the Undecidability of Freeness of Matrix Semigroups},
\newblock Int. J. Algebra Comput. 9 (3--4), 295--305 (1999)

\bibitem{CDM}
A.~S.\,Cavaretta,  W.\,Dahmen, and C.A.\,Micchelli,
\newblock  {\em Stationary
subdivision},
\newblock Mem. Amer. Math. Soc. 93 (1991), no.~453.

\bibitem{Cha}
G.M.\,Chaikin,
\newblock {\em An algorithm for high speed curve generation},
\newblock Comp. Graphics and Image. Proc., 3 (1974), 346--349.

\bibitem{dR1}
G.\,de Rham,
\newblock {\em Sur une courbe plane},
\newblock  J. Math. Pures Appl. 35 (1956), 25–-42.

\bibitem{dR2}
G.\,de Rham,
\newblock {\em Sur les courbes limites de polygones obtenus par trisection},
\newblock  Enseignement Math. 5 (1959), 29--43.

\bibitem{D}
G.A.\,Derfel,
\newblock {\em Probabilistic method for a class of functional-differential equations}
\newblock Ukrain. Math. J., 41 (1989), 1137-–1141.

\bibitem{DGL}
N.\,Dyn, J.A.\,Gregory, and D.\,Levin,
\newblock {\em Analysis of uniform binary
subdivision schemes for curve design},
\newblock Constr. Approx. 7 (1991), 127--147.

\bibitem{DL}
N.\,Dyn and D.\,Levin,
\newblock {\em Subdivision schemes in geometric modelling},
\newblock Acta Numer., 11 (2002), 73--144.

\bibitem{ErdGra}
P.\,Erdos, R.\,Graham,
\newblock {\em Old and new problems and results in combinatorial number theory},
\newblock L'Enseignement Math., Monogr., Vol. 28, 1980.

\bibitem{FS}
D.-J.\,Feng and N.\,Sidorov,
\newblock {\em Growth rate for beta-expansions},
\newblock Monatsh. Math.
162 (2011), 41–-60

\bibitem{Hon1}
J.\,Honkala,
\newblock {\em Unique representation in number systems and L codes},
\newblock Discrete Applied Mathematics 4 (1982), 229--232.

\bibitem{Hon2}
J.\,Honkala,
\newblock {\em On Number Systems with Finite Degree of Ambiguity},
\newblock Information and Computation 145 (1998), 51--63.

\bibitem{H}
J.\,E.\,Hutchinson,
\newblock {\em  Fractals and self-similarity},
\newblock Indiana Univ. Math.~J., 30 (1981), no 5, 713--747.

\bibitem{JT}
J.\,Jankauskas and J.M.\,Thuswaldner, 
\newblock {\em  Rational matrix digit systems}, 
\newblock Linear and Multilinear Alg.,, 71 (2023), 
no 10,  1606–-1639. 

\bibitem{K81}
D.~A.\,Klarner,
\newblock {\em An Algorithm to Determine When Certain Sets Have O-Density},
\newblock Journal of Algorithms 2 (1981), 31--43.

\bibitem{K82}
D.A.\,Klarner,
\newblock {\em A Sufficient Condition for Certain Semigroups to Be Free''}
\newblock J. Algebra 74, 140--148 (1982).

\bibitem{KBS}
D.A.\,Klarner, J.C.\,Birget, W.\,Satterfield,
\newblock {\em On the undecidability of the freeness of integer matrix semigroups},
\newblock International J.\ of Algebra and Computation 1 (1991) 223-226.

\bibitem{KT}
A.\,Kolpakov, A.\,Talambutsa,
\newblock {\em On free semigroups of affine maps on the real line},
\newblock Proc. Amer. Math. Soc., 150:6 (2022), 2301–2307.

\bibitem{Lag}
J.\,C.~Lagarias,
\newblock {\em ``Erd\"os, Klarner, and the $3x+1$ problem''},
\newblock Amer. Math. Monthly 123 (8), 753--776 (2016)

\bibitem{LW}
J.\,C.~Lagarias and Y.Wang,  
\newblock {\em  Integral self-affine tiles in $\R^n$. I. Standard and nonstandard digit sets} 
\newblock J. London Math. Soc., 54 (1996), no 1, 161-–179.

\bibitem{Li}
J.-L.\, Li,  
\newblock {\em  Digit sets of integral self-affine tiles with prime determinant}, 
\newblock Studia Mathematica, 
177 (2006), no 2, 183--194.  

\bibitem{MSW}
H.A.\,Maurer, A.\,Salomaa, D.\,Wood,
\newblock{\em L codes and number systems},
\newblock Theoretical Computer Science 22 (1983), 331--346.

\bibitem{MP}
C.A.\,Micchelli and H.\,Prautzsch,
\newblock {\em Uniform refinement of curves},
\newblock  Lin. Alg. Appl. 114/115 (1989),  841--870.

\bibitem{NPS}
I.Ya.\,Novikov, V.Yu.\,Protasov, and M.A.\,Skopina,
\newblock {\em Wavelets theory},
\newblock AMS, Translations Mathematical
Monographs, 239 (2011).

\bibitem{P00b}
V.\,Yu.~Protasov,
\newblock {\em Asymptotic behaviour of the partition function},
\newblock Sb. Math., 191 (2000), No 3--4, 381--414.

\bibitem{P00p}
V.\,Yu.~Protasov,
\newblock {\em Refinement equations with nonnegative coefficients},
\newblock J. Fourier Anal. Appl., 6 (2000), No 1,  55--78.

\bibitem{P04b}
V.\,Yu.~Protasov,
\newblock {\em On the asymptotics of the binary
partition function}, Math. Notes,  76 (2004), No 1,
151--156.

\bibitem{P17}
V.\,Yu.~Protasov, \newblock{\em The Euler binary partition function and subdivision schemes}, Mathematics of Computation, 86 (2017), No 305, 1499--1524.

\bibitem{Z}
O.K.\,~Zakusilo,
\newblock {\em Some properties of classes $L_c$ of limit distribution},
\newblock  Teoria Veroyatnosti i Mat. Statistika, 15 (1976), 68--73.

\end{thebibliography}
\end{document}